\documentclass[11pt]{article}
\usepackage[top=1.5in, bottom=1in, left=1in, right=1in]{geometry}
 \usepackage[backref,colorlinks,linkcolor=red,anchorcolor=green,citecolor=blue]{hyperref}
\usepackage{amsfonts,amssymb}
\usepackage{amsmath}
\usepackage{graphicx}
\usepackage{cite}
\usepackage{enumerate}
\usepackage{caption}
\usepackage{subcaption}
\usepackage{booktabs} 
\usepackage{multirow}
\usepackage[ruled]{algorithm2e}

\newcommand{\I}{\mathcal{I}}
\newcommand{\J}{\mathcal{J}}
\newcommand{\Ir}{\widetilde{\I}}
\newcommand{\Is}{\hat{\I}}

\newcommand{\B}{\mathcal{B}}
\newcommand{\R}{\mathcal{R}}

\newcommand{\A}{{A}}

\newcommand{\AB}{\Bar{\A}}
\newcommand{\ABB}{\widetilde{\A}}
\newcommand{\E}{{M}}
\newcommand{\F}{{F}}

\newcommand{\K}{{K}}
\renewcommand{\L}{{L}}
\newcommand{\T}{{T}}
\newcommand{\Z}{{Z}}
\newcommand{\C}{{C}}

\allowdisplaybreaks

\newcommand{\rev}[1]{#1} 

\allowdisplaybreaks
\begin{document}

\date{}

\title{Hierarchical Interpolative Factorization Preconditioner \\for Parabolic Equations}


\author{Jordi Feliu-Fab\`{a}\thanks{Institute for Computational and Mathematical Engineering,
    Stanford University, Stanford, CA 94305, (jfeliu@stanford.edu).}  \and Lexing
  Ying\thanks{Department of Mathematics and Institute for Computational and Mathematical
    Engineering, Stanford University, Stanford, CA 94305, (lexing@stanford.edu)}}

\pagestyle{myheadings} \markright{\hfill HIF preconditioner for parabolic equations \hfill}

\maketitle

\begin{abstract}
  This note proposes an efficient preconditioner for solving linear and semi-linear parabolic
  equations.  With the Crank-Nicholson time stepping method, the algebraic system of equations at
  each time step is solved with the conjugate gradient method, preconditioned with hierarchical
  interpolative factorization. Stiffness matrices arising in the discretization of parabolic
  equations typically have large condition numbers, and therefore preconditioning becomes essential,
  especially for large time steps. We propose to use the hierarchical interpolative factorization as
  the preconditioning for the conjugate gradient iteration.  Computed only once, the hierarchical
  interpolative factorization offers an efficient and accurate approximate inverse of the linear
  system. As a result, the preconditioned conjugate gradient iteration converges in a small number
  of iterations. Compared to other classical exact and approximate factorizations such as Cholesky
  or incomplete Cholesky, the hierarchical interpolative factorization can be computed in linear
  time and the application of its inverse has linear complexity. Numerical experiments demonstrate
  the performance of the method and the reduction of conjugate gradient iterations.
\end{abstract}


\section{Introduction}\label{introduction}

This note is concerned with the numerical solution of parabolic equations of the form
\begin{align}
  \label{eq1}
  \frac{\partial u(x,t)}{\partial t} = \nabla\cdot\big(a(x)\nabla u(x)\big)+r\big(u(x,t)\big),\quad
  x\in\Omega\subset\mathbb{R}^{d},
\end{align}
in two and three dimensions, with appropriate boundary conditions on $\partial\Omega$ and initial
conditions $u(x,0) = u_0(x)$. Here, $a(x)>0$ is the coefficient field of the diffusion operator and
$r\big(u(x,t)\big)$ is the reaction term. We are interested on approximating the unknown field
$u(x,t)$. Such reaction-diffusion equations can model a great variety of physical phenomena, such as
heat conduction with internal heat generation, population dynamics \cite{fisher,dynamics} or pattern
formation \cite{turing} in biology.

The most common spatial discretizations for solving \eqref{eq1} are finite difference and finite
element methods. Such a spatial discretization results in a time-dependent system of form
\begin{align}
  \label{eq2}
  \frac{\partial u(t)}{\partial t} = M u(t) + r(t),
\end{align}
where $u(t)\in\mathbb{R}^{N}$ and $r(t)\in\mathbb{R}^{N}$ are the spatial discretizations of
$u(x,t)$ and $r\big(u(x,t)\big)$ at time $t$, respectively, and $N$ is the number of degrees of
freedom (DOFs) in the spatial discretization. The {\em stiffness} matrix $M\in \mathbb{R}^{N\times
  N}$ is the discretization of the diffusion term.

The classical approach to solve \eqref{eq2} consists of discretizing in time and evolving the
numerical solution at successive time steps using a time marching method. For instance, an explicit
scheme to approximate the solution at every time step successively can be obtained by using the
forward Euler method in time,
\[
u^{k+1} = (I + M \Delta t) u^k + \Delta t \cdot r^k.
\]
However, such an approach would require selecting a very small time step $\Delta t$ in order to
satisfy the stability condition $\Delta t \leq \frac{1}{2}\Delta x^2\max_{x\in\Omega} a(x)$, where
$\Delta x$ is the spacing of the spatial grid.

Alternatively, one can overcome the stability condition and use larger time steps by using an
implicit scheme such as the first-order backward Euler or the second-order Crank-Nicolson methods,
which are unconditionally stable. However, such methods introduce another challenge by requiring to
solve a system of equations at each time step. This note focuses on the Crank-Nicolson method
\begin{align}
  \label{eq_CN}
\left(I-\frac{\Delta t}{2}M\right) u^{k+1} = \left(I+\frac{\Delta t}{2}M\right) u^k + \Delta t\cdot r^k.
\end{align}
There are several ways to solve the system of equations. The most direct way is to use exact
factorization methods, such as Cholesky decomposition, which would be prohibitively expensive for
large problem sizes. Alternatively, one can use iterative methods such as conjugate gradient (CG)
which takes $O(nnz(M))=O(N)$ cost per iterations. However, since the matrix $I-\frac{\Delta t}{2}M$
is typically ill-conditioned, the number of iterations can be quite large, thus requiring the use of
a preconditioner. Finding a good preconditioner is itself a challenging task that has been vastly
studied in the literature. Among others, we highlight here the incomplete Cholesky factorization
\cite{ichol} and multifrontal method based algorithms coupled with hierarchical matrices
\cite{Xia:2009} and skeletonization \cite{HIFDE,PHIF,cambier,Klockiewicz2019SparseHP,Abey,Abey2}.

The main goal of this note is to describe a new efficient preconditioner based on a version of the
{\em hierarchical interpolative factorization} (HIF) described in \cite{PHIF} to reduce the number
of CG iterations at each time step. Due to the accuracy and efficiency of this version of HIF, the
preconditioned CG iteration at each time step converges in a small number of iterations. We
demonstrate the effectiveness of this approach of solving \eqref{eq1} by studying several 2d and 3d
numerical examples.

\paragraph{Other approaches.}\label{sec:background}
In recent years, several exponential integrator based models for time discretization have been
proposed to solve parabolic equations, based on the integration factor method \cite{Lili} and the
exponential time differencing method \cite{HUANG2019257,zhu}.

Boundary integral formulations provide an alternative to the PDE-based approach.  In the case where
the diffusion coefficient $a(x)$ is constant, one can make use of parabolic potential theory
\cite{TAUSCH2007956, Wang2018Gauss} to obtain a boundary integral equation including time convolution of the form
\begin{align}
  \label{eq_ie}
  u(x,t) = \int_{\Omega(0)}{G(x-y,t)u_0(y)dy} +
  \int_{0}^{t}\int_{\Omega(\tau)}{G(x-y,t-\tau)r(y,\tau)dyd\tau},
\end{align}
where $G(x,t)$ is the free-space Green's function for the heat equation in $d$ dimensions
\begin{align}
  \label{eq_heatkernel}
  G(x,t)=\frac{e^{-\|x\|^2/4t}}{(4\pi t)^{d/2}}.
\end{align}
One can then use fast algorithms to evaluate layer potentials, such as the hierarchical
interpolative factorization \cite{HIFIE}, the fast multipole method \cite{GREENGARD1987325} as done
by Messner et al. \cite{MESSNER201415} or the fast Gauss transform
\cite{FG,Tausch2009MultidimensionalFG} used by Wang et al. \cite{Wang2019FastIE} in this context.
Other works include \cite{barnett,Wang2019_2}. This approach has many advantages, such as time
stability, reduction of DOFs to the boundary and the availability of fast linear algorithms for the
evaluation of the layer potentials. However, so far, the classical potential theory is restricted to
the set of reaction-diffusion equations with constant coefficients.

\section{Algorithm}\label{algorithm}

This section reviews the version of HIF introduced in \cite{PHIF} and then describes the new
preconditioner. Throughout the note, the following notation are used: the uppercase letters ($A$,
$F$, $M$, etc.)\ denote matrices; the calligraphic letters ($\I$, $\B$, $\R$, etc.)\ denote sets of
indices, associated to DOFs; $A_{\I\B}$ refers to the restriction of $A$ to the $|\I|\times|\J|$
submatrix with rows indexed by $\I$ and columns indexed by $\B$; the notation $\{\I_{i}\}_{i=1}^{p}$
represents a collection of $p$ disjoint sets of indices $\I_{i}$ for $i=1,\dots, p$.

Consider the PDE \eqref{eq1} on $\Omega=(0,1)^2$ with the Dirichlet boundary conditions and initial
condition $u(x,0)=u_0(x)$. We perform finite difference dicretization via the standard five-point
stencil over a uniform grid with step size $h$, resulting in $N=(n-1)^2$ DOFs, with $n=1/h=2^Lm$ for
some integers $L$ and $m$. Each DOF corresponds to the solution $u_j = u(x_j)$ at grid points
$x_j=(hj_1,hj_2)$, with $j=(j_1,j_2)$ and $1<j_1,j_2<n-1$. The resulting matrix $M$ corresponding to
the discretization of the diffusion term $\nabla a(x)\nabla$ in \eqref{eq1} is sparse and symmetric
negative definite.

\subsection{Crank-Nicolson scheme}
\label{CrankNicolson}

Since the diffusive term is the leading term, we use the Crank-Nicolson scheme for the diffusive
term and an explicit scheme for the reaction term, with time step $\Delta t$. This leads to the
algebraic system of equations \eqref{eq_CN}, which can be solved successively to evolve the
numerical solution at future time steps,
\begin{align}
  \label{eq_CN2}
  u^{k+1} =\left(I-\frac{\Delta t}{2}M\right)^{-1} \left[\left(I+\frac{\Delta t}{2}M\right)u^k + \Delta t \cdot r^k\right].
\end{align}
Since the unconditionally stable Crank-Nicolson scheme is second order in time and space, it allows
for larger time steps than Backward or Forward Euler methods. Since $A = I-\frac{\Delta t}{2} M$ is
sparse and symmetric-positive-definite (SPD) one can use conjugate gradient to solve \eqref{eq_CN2}
as opposed to classical direct solvers which are more expensive. However, the number of iterations
scale with $\kappa(A)$, the condition number of $A$. For small time steps, $\kappa(A)$ may be close
to 1. However, for the large time steps adopted here, $\kappa(A)$ is comparable to the condition
number of $M$, which is quite large as $M$ is ill-conditioned. Preconditioning then becomes an
essential step in order to reduce the number of CG iterations. The HIF has been shown to
significantly reduce the number of CG iterations for the Poisson equation with variable coefficient
\cite{PHIF} and we apply it to precondition $A = I-\frac{\Delta t}{2} M$.

\subsection{Recursively preconditioned hierarchical interpolative factorization}\label{PHIF}

Here we briefly review the {\em recursively preconditioned hierarchical interpolative factorization}
(PHIF) proposed in \cite{PHIF}, which is an improved version of the original HIF described in
\cite{HIFDE}. Given a sparse SPD matrix $A$, PHIF produces a fast factorization that can be seen as
a multilevel generalized Cholesky decomposition, which alternates among block Gaussian elimination,
block Jacobi preconditioning and skeletonization at different levels.

 One starts by defining a uniform quad-tree with $L$ levels,  that partitions the domain $\Omega$ into
$p^{\ell}=2^{L-\ell}\times2^{L-\ell}$ square cells at each level $\ell=0,1,\hdots,L$. The leaves
correspond to level $\ell=0$ and the root to level $\ell=L$. Throughout the factorization we denote
\textit{active} DOFs the DOFs that have not yet been decoupled. We set $\A_0 = \A$ and start
factorizing the matrix by decoupling DOFs starting from the leaves level $\ell=0$ up to level $L-1$
with the following three steps at each level $\ell$:
\vspace{0.5cm}

\begin{enumerate}[(1)]

\item \textbf{Cell elimination.} For each cell $1<i<p_\ell$ at level $\ell$ we decouple the interior
  active DOFs indexed by $\I_{\ell,i}$ as follows.  Up to a permutation matrix $A_{\ell}$ can be
  written as
\begin{align}
\A_{\ell} = \begin{bmatrix}
  \A_{\I_{\ell,i}\I_{\ell,i}}&\A_{\B\I_{\ell,i}}^T&\\
  \A_{\B\I_{\ell,i}}&\A_{\B\B}&\A_{\R\B}^T\\
  &\A_{\R\B}&\A_{\R\R}\\
\end{bmatrix},
\label{eqn:matrix-3blocks}
\end{align}
where $\I_{\ell,i}$ represents the DOFs inside that cell, $\B$ the
DOFs on the boundary (i.e. edges and corners of the cell), and $\R$ the remaining DOFs. For clarity purposes, hereafter we use sub-index $\ell$ to denote matrices at level $\ell$ and drop the sub-index $\ell$ when referring to a particular block of the matrix, for instance we use $A_{\B\B}$ to refer to $A_{\ell_{\B\B}}$.   Let $A_{\I_{\ell,i},\I_{\ell,i}} =
L_{\I_{\ell,i}}L_{\I_{\ell,i}}^T$, Gaussian elimination leads to
\begin{align} \label{eq2_3}
\E_{\I_{\ell,i}}^{T}\A_{\ell} \E_{\I_{\ell,i}}=\begin{bmatrix}
  I&&\\
  &X_{\B\B}&\A_{\R\B}^T\\
  &\A_{\R\B}&\A_{\R\R}\\
\end{bmatrix}, \quad \E_{\I_{\ell,i}}=
  \begin{bmatrix}
    \L_{\I_{\ell,i}}^{-T} \hspace*{0.1cm}&\hspace*{0.1cm} -A_{\I_{\ell,i}\I_{\ell,i}}^{-1} A_{\B\I_{\ell,i}}^{T}&\hspace*{0.1cm}\\
    &I\\
    &&\hspace*{0.1cm}I
  \end{bmatrix},
\end{align}
where $X_{\B\B}=\A_{\B\B}-\A_{\B\I_{\ell,i}}\A_{\I_{\ell,i}\I_{\ell,i}}^{-1}\A_{\B\I_{\ell,i}}^T$.

Performing block eliminations for each cell $1<i<p_\ell$ at level $\ell$ of the quadtree leads to
\begin{align}\label{eq3_1}
  \AB_{\ell}=\E_{\ell}^{T}\A_{\ell}\E_{\ell}, \quad \E_{\ell}=\prod_{i=1}^{p_\ell} \E_{\I_{\ell,i}}.
 \end{align}
After block elimination with respect the collection of index sets $\{ \I_{\ell,i} \}_{i =
  1}^{p_\ell}$, the DOFs inside the cells have been decoupled from those in the edges and corners at
level $\ell$. Therefore, the remaining active DOFs are those in the edges and corners and we only need to
continue factorizing the matrix $\AB_\ell$ restricted on these active DOFs.

\item \textbf{Block Jacobi preconditioning.}  Let $r_\ell$ be the number of edges and corners at
  level $\ell$, and $\{ \I_{\ell, i} \}_{i = 1}^{r_\ell}$ be the collection of corresponding index
  sets for the active DOFs. For a given edge or corner $i$, up to a permutation, $\AB_\ell$ can be
  written as
\begin{align}
\AB_\ell = \begin{bmatrix}
  \AB_{\I_{\ell,i}\I_{\ell,i}}&\AB_{\B\I_{\ell,i}}^T&\\
  \AB_{\B\I_{\ell,i}}&\AB_{\B\B}&\AB_{\R\B}^T\\
  &\AB_{\R\B}&\AB_{\R\R}\\
\end{bmatrix},
\label{eqn:matrix-3block}
\end{align}
where $\I_{\ell,i}$ represents the DOFs in the edge/corner $i$,
$\B$ the DOFs on the edges and corners connected to edge/corner $i$ in $\AB_{\ell}$, and $\R$ the remaining active DOFs.
Then a \emph{rescaling} of edge/corner $i$ can be performed using the Cholesky decomposition
$\AB_{\I_{\ell,i}\I_{\ell,i}} = L_{\I_{\ell,i}} L_{\I_{\ell,i}}^T$ as
\begin{align}
 \C_{\I_{\ell,i}}^T \AB_{\ell} \C_{\I_{\ell,i}} =
 \begin{bmatrix}
  I &  L_{\I_{\ell,i}}^{-1}\AB_{\B\I_{\ell,i}}^T\\
  \AB_{\B\I_{\ell,i}} L_{\I_{\ell,i}}^{-T} & \AB_{\B\B} & \AB_{\R\B}^T\\
  & \AB_{\R\B} & \AB_{\R\R}
 \end{bmatrix}, \quad \C_{\I_{\ell,i}} =
 \begin{bmatrix}
  L_{\I_{\ell,i}}^{-T}\\
  & I\\
  && I
 \end{bmatrix} \in \mathbb{R}^{N \times N},
 \label{eqn:jacobi-precond}
\end{align}

If we perform this preconditioning for each edge and corner in level $\ell$, we obtain a block Jacobi
preconditioning that yields
\begin{align} \label{eq3_2}
  \ABB_{\ell} = \C_{\ell}^T \AB_{\ell} \C_{\ell}, \quad \C_{\ell} = \prod_{i = 1}^{r_\ell} \C_{\I_{\ell, i}},
\end{align}
where $\C_{\ell}$ is a block diagonal matrix, up to a permutation, since $\{ \I_{\ell, i} \}_{i =
  1}^{r_\ell}$ is a collection of disjoint index sets.  The diagonal blocks of the resulting matrix
$\ABB_{\ell}$ are identity matrices and the set of active DOFs remains unchanged.
 
\item \textbf{Edge skeletonization.}  Let $q_\ell$ be the number of edges at level $\ell$, and
  $\{\I_{\ell, i} \}_{i = 1}^{q_\ell}$ the collection of corresponding index sets.

For a given edge $i$ with DOFs indexed by $\I_{\ell, i}$ skeletonization is performed as
follows. For clarity we drop the subindex in $\I_{\ell, i}$ in the remainder of the explanation
for the edge skeletonization step.  Up to a permutation one can write
\begin{align}\label{eq2_4}
\ABB_{\ell} =
\begin{bmatrix}
  I&\ABB_{\R\I}^T\\ \ABB_{\R\I}&\ABB_{\R\R}\\
\end{bmatrix}
\end{align}
where $I$ are the DOFs in edge $i$ at level $\ell$ and $R$ are the rest of active DOFs. Assume $\ABB_{\R\I}
\in\mathbb{R}^{N_\R \times N_\I}$ has numerical rank $k$ to relative precision
$\epsilon$. Then, we can use interpolative decomposition (ID) \cite{ID}, to get a partition of
$\I=\Ir\cup\Is$ into \emph{skeleton} $\Is$ and \emph{redundant} $\Ir$ DOFs, and approximate the
redundant columns of $\ABB_{\R\I}$ by a linear combination of its skeleton columns such that
\begin{align}\label{eq2_5}
\ABB_{\R\Ir} = \ABB_{\R\Is}\T_{\I} + E_{\I}, \quad \|E_{\I}\| = O(\epsilon\|\ABB_{\R\I}\|),
\end{align}
with $\T_{\I}\in\mathbb{R}^{k\times (N_\I-k)}$. 
Up to a permutation, using ID we can approximately zero out the redundant columns leading to
\begin{align}
  \label{eq_T}
  \Z_{\I}^{T}\ABB_{\ell}\Z_{\I}\approx\left[{\begin{array}{cc|c}
        I+\T_\I^{T}\T_\I\hspace*{0.2cm}&-\T_\I^{T}&\\
        \vspace*{0.1cm} -\T_\I&I&\ABB_{\R\Is}^T\\ \hline
       \\[-2.5ex] &\ABB_{\R\Is}&\ABB_{\R\R}\\
    \end{array}}\right],
  \quad
  \Z_{\I}=\left[{\begin{array}{cc|c}
        I&&\\
        -\T_{\I}&I&\\ \hline
        &&I\\
    \end{array}}\right]\in\mathbb{R}^{N\times N},
\end{align}

Now, we can decouple redundant DOFs $\Ir$ in the edge by performing Gaussian elimination using the
Cholesky decomposition $I+\T_\I^{T}\T_\I=L_{\Is}L_{\Is}^{T}$,
\begin{align}\label{eq2_9}
\E_{\Ir}^{T}\Z_{\I}^{T}\ABB_{\ell}\Z_{\I}\E_{\Ir}\approx\left[{\begin{array}{cc|c}
      I&&\\
      &X_{\Is\Is}\hspace*{0.1cm}&\hspace*{0.1cm} \ABB_{\R\Is}^T \vspace*{0.1cm}\\  \hline
       \\[-2.5ex]
      &\ABB_{\R\Is}\hspace*{0.1cm}&\hspace*{0.1cm}\ABB_{\R\R}\\
  \end{array}}\right], \quad \E_{\I}=
  \begin{bmatrix}
    L_{\Is}^{-T} &\hspace*{0.1cm} (I+\T_\I^{T}\T_\I)^{-1} T_{\I}^{T}&\hspace*{0.2cm}\\
    &I&\\
    &&\hspace*{0.2cm}I
  \end{bmatrix} 
\end{align}
with $X_{\Is\Is}=I-T_{\I}(I+\T_\I^{T}\T_\I)^{-1}T_{\I}^{T}$.

Performing skeletonization for each edge $i$ with set of indices $\I_{\ell,i}$ gives
\begin{align} \label{eq3_3}
  \A_{\ell+1}\approx\K_{\ell}^{T}\ABB_{\ell}\K_{\ell}, \qquad
  \K_{\ell}=\prod_{i=1}^{q_\ell}\K_{\I_{\ell,i}}, \quad \K_{\I_{\ell,i}}
  =\Z_{\I_{\ell,i}}\E_{\I_{\ell,i}}.
\end{align}

The remaining active DOFs are now the skeleton DOFs in the edges and the corners at level $\ell$. We can now move
to the following level on the tree and perform the same three steps.
\end{enumerate}

\vspace{0.5cm}

After performing these three steps for all levels $\ell=0,\hdots,L-1$, i.e. once we are in the root
level of the tree, the resulting matrix is
\begin{align} \label{eq3_4}
 \A_L \approx  R_{L-1}^T \cdots R_0^T \A R_0 \cdots R_{L-1}, \quad R_\ell = \E_\ell \C_\ell \K_\ell
\end{align}
which is everywhere the identity except in the block indexed by the remaining active DOFs at the
root level. As opposed to the nested dissection multifrontal factorization, the frontal matrix at the
root is small since sparsification on the separator fronts has been performed throughout all the levels.  Now, we can
approximate the original matrix as
\begin{align} 
 \A \approx \F = R_{0}^{-T} \cdots R_{L-1}^{-T}\A_{L}R_{L-1}^{-1} \cdots R_{0}^{-1},
 \label{eqA}
\end{align}
and its inverse as
\begin{align}
 \A^{-1} \approx \F^{-1} = R_{0} \cdots R_{L-1}\A_{L}^{-1}R_{L-1}^{T} \cdots R_{0}^{T}.
 \label{eqAinv}
\end{align}

The factors $R_\ell=\E_\ell \C_\ell \K_\ell$ are easily invertible since $\E_\ell, \C_\ell, \K_\ell$
are block diagonal up to a permutation, with each block being triangular. Therefore, the
factorization can be viewed as a generalized Cholesky decomposition.

In the 3d case, the algorithm for PHIF is similar to the 2d case, except for a few details. Instead
of square cells, the domain is partitioned into cube cells. At each level $\ell$, Block Jacobi
preconditioning is performed in faces, edges and corners, while skeletonization is performed on
faces (although it can also be performed in edges to lower the computational complexity of the
factorization).

\begin{figure}
\centering \captionsetup[subfigure]{labelformat=empty}
\begin{subfigure}{0.225\textwidth}
 \centering
 \includegraphics[width=\textwidth,trim=4.30cm 2.05cm 3.35cm 2.0cm,clip]{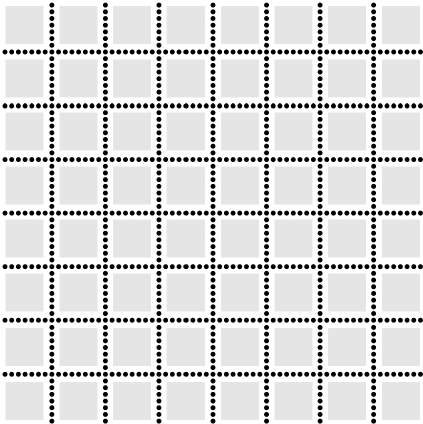}
 \caption{$\ell=0$ cell elimination}
\end{subfigure}
\begin{subfigure}{0.225\textwidth}
 \centering
 \includegraphics[width=\textwidth,trim=4.30cm 2.05cm 3.35cm 2.0cm,clip]{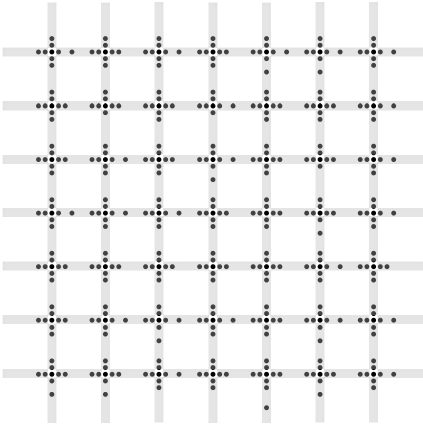}
 \caption{$\ell=0$ skeletonization}
\end{subfigure}
\begin{subfigure}{0.225\textwidth}
 \centering
 \includegraphics[width=\textwidth,trim=4.30cm 2.05cm 3.35cm 2.0cm,clip]{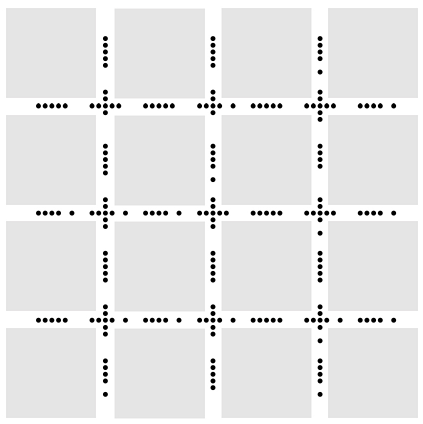}
 \caption{$\ell=1$ cell elimination}
\end{subfigure}\\
\vspace{2mm}
\begin{subfigure}{0.225\textwidth}
 \centering
 \includegraphics[width=\textwidth,trim=4.30cm 2.05cm 3.35cm 2.0cm,clip]{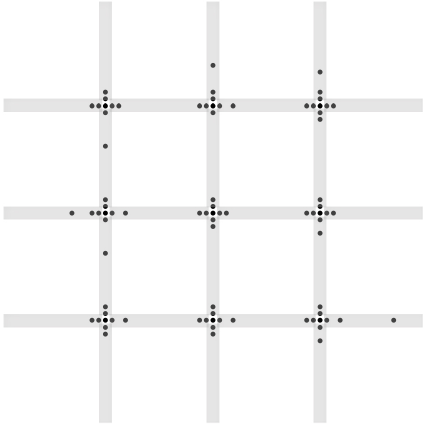}
 \caption{$\ell=1$ skeletonization}
\end{subfigure}
\begin{subfigure}{0.225\textwidth}
 \centering
 \includegraphics[width=\textwidth,trim=4.30cm 2.05cm 3.35cm 2.0cm,clip]{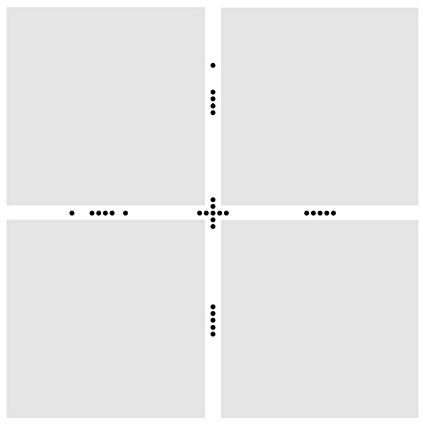}
 \caption{$\ell=2$ cell elimination}
\end{subfigure}
\begin{subfigure}{0.225\textwidth}
 \centering
 \includegraphics[width=\textwidth,trim=4.30cm 2.05cm 3.35cm 2.0cm,clip]{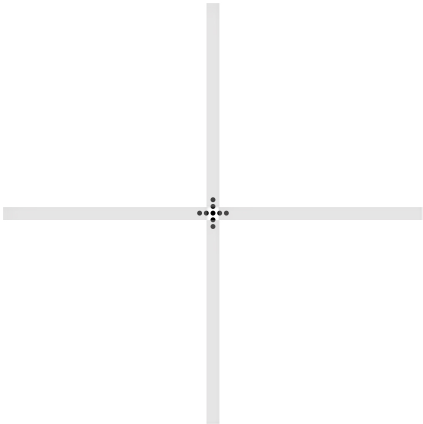}
 \caption{$\ell=2$ skeletonization}
\end{subfigure}
\caption{Active DOFs at each level $\ell$ of PHIF in 2d, depicted with black dots. The square cells
  and edges at each level, from which DOFs have been eliminated, are represented in gray for cell
  elimination and skeletonization steps respectively.}
\label{mHIFDE-fig}
\end{figure}

\begin{figure}
\centering \captionsetup[subfigure]{labelformat=empty}
\begin{subfigure}{0.225\textwidth}
 \centering
 \includegraphics[width=\textwidth,trim=5.0cm 3.1cm 4.5cm 2.3cm,clip]{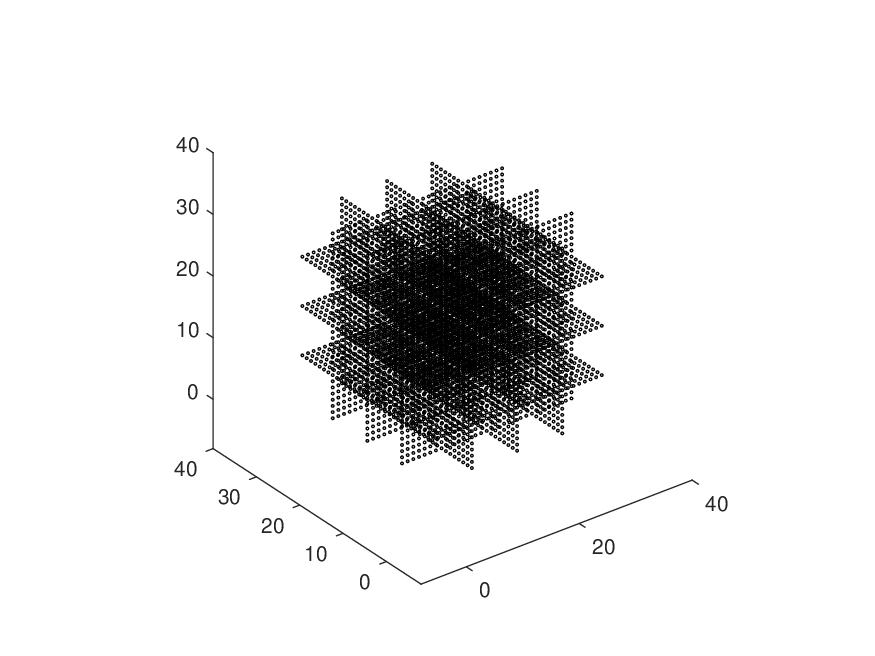}
 \caption{$\ell=0$ cell elimination}
\end{subfigure}
\begin{subfigure}{0.225\textwidth}
\centering
\includegraphics[width=\textwidth,trim=5.0cm 3.1cm 4.5cm 2.3cm,clip]{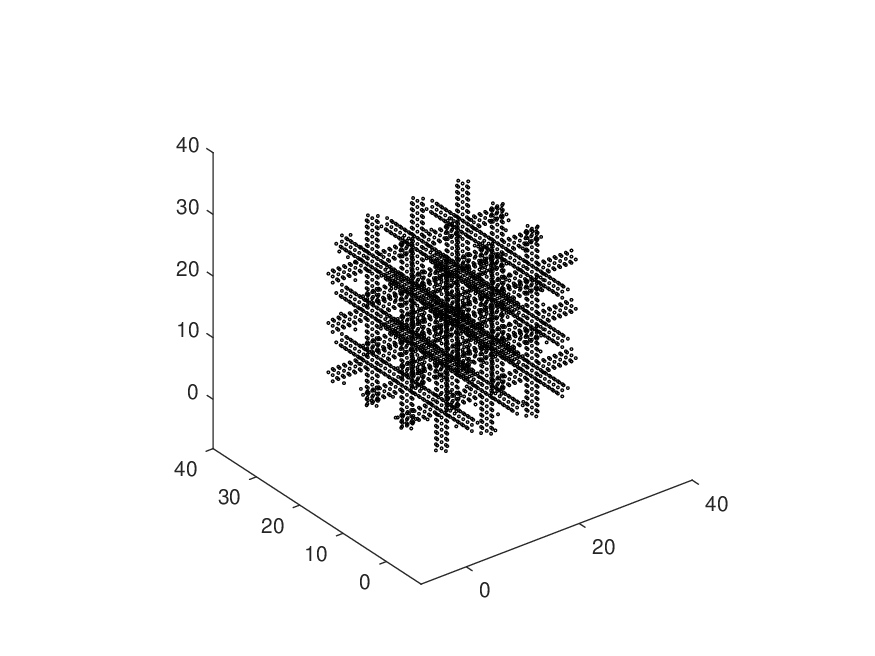}
\caption{$\ell=0$ skeletonization}
\end{subfigure}
\begin{subfigure}{0.225\textwidth}
 \centering
 \includegraphics[width=\textwidth,trim=5.0cm 3.1cm 4.5cm 2.3cm,clip]{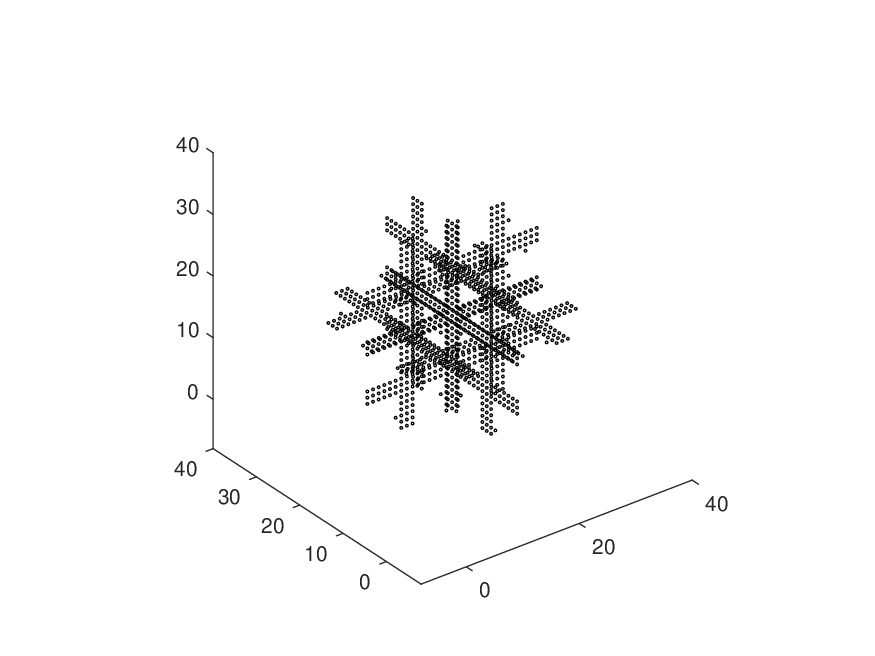}
 \caption{$\ell=1$ cell elimination}
\end{subfigure}
\begin{subfigure}{0.225\textwidth}
 \centering
 \includegraphics[width=\textwidth,trim=5.0cm 3.1cm 4.5cm 2.3cm,clip]{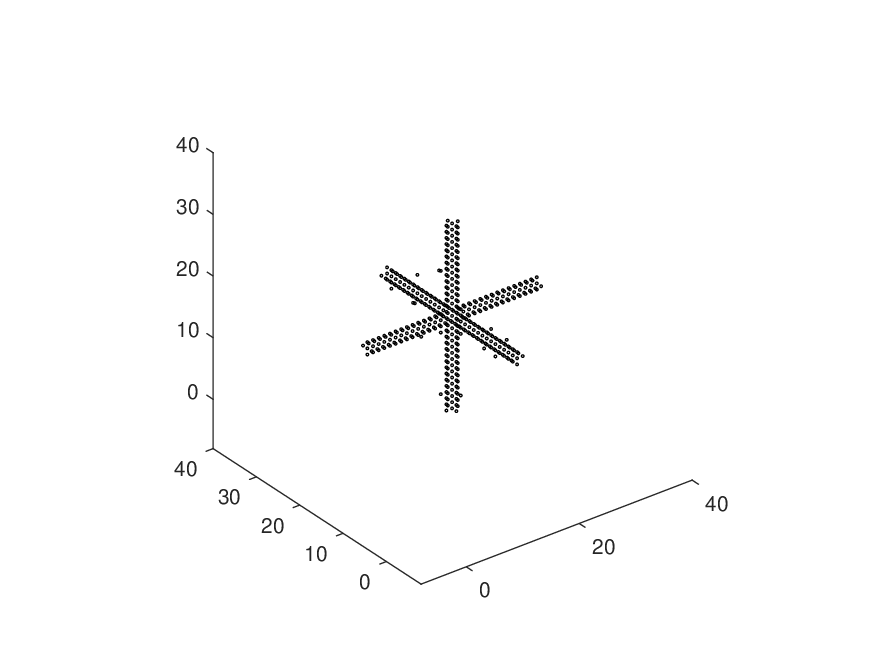}
 \caption{$\ell=1$ skeletonization}
\end{subfigure}\\
\vspace{3mm}
\begin{subfigure}{0.2\textwidth}
 \centering
 \includegraphics[width=\textwidth,clip]{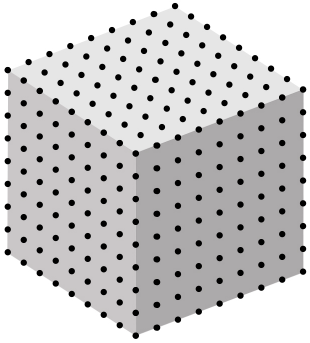}
 \caption{DOFs on faces}
\end{subfigure}
\qquad
\begin{subfigure}{0.2\textwidth}
 \centering
 \includegraphics[width=\textwidth,clip]{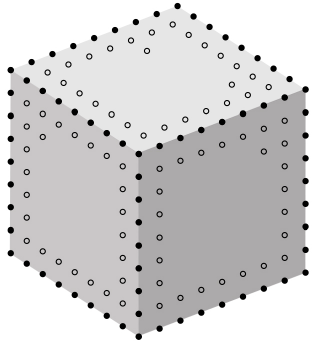}
 \caption{Skeletons}
\end{subfigure}
\caption{Top: Active DOFs at each level $\ell$ of PHIF in 3d. Bottom: Detail of DOFs on the faces of a cubic cell before and after skeletonization, with skeleton DOFs represented by hollow circles.}
\label{mHIFDE-fig3d}
\end{figure}

See Figures \ref{mHIFDE-fig} and \ref{mHIFDE-fig3d} for an illustration of the active DOFs at
different steps and levels of the factorization process in 2d and 3d respectively. The Block Jacobi
preconditioning step has been omitted in the illustrations because it doesn't change the active
DOFs.

\subsection{Crank-Nicolson scheme with PHIF preconditioned CG}\label{alg_CN}

We now proceed to describe the algorithm for evolving the numerical solution of \ref{eq1}. First, we
compute the PHIF factorization $F \approx A = (I-\frac{\Delta t}{2}M)$. Since $A$ is a sparse SPD
matrix, the PHIF can be computed in linear time to obtain $F^{-1}\approx A^{-1}$, and matrix-vector
multiplications with $F^{-1}$ require only $O(N)$ work.

In each time step of the Crank-Nicolson method, the linear system of equation \eqref{eq_CN} is
solved by CG, with $F^{-1}$ in its factorized form \eqref{eqAinv} used as the preconditioner. The
algorithm is presented in Algorithm \ref{alg}, where $\text{pcg}(A,b,P)$ represents the
preconditioned conjugate gradient with preconditioner $P$ that attempts to solve for $Ax=b$.

\begin{algorithm}[htb]
  Construct PHIF $F$ for $A = (I-\frac{\Delta t}{2}M)$ with some tolerance $\epsilon$\;
  Initialize $u^0=u_0$\;
  \For{$k = 1: n_{steps}$}{
    $g = (I+\frac{\Delta t}{2}M)u^{k-1} + \Delta t \cdot r^{k-1}$\; 
    $u^k = \text{pcg}(A,g, F^{-1})$ (see \eqref{eqAinv} for the factorization of $F^{-1}$)
  }
  \caption{Numerical solution of reaction-diffusion equations}
  \label{alg}
\end{algorithm}

\section{Numerical Results}\label{Numerical}

In this section, we demonstrate the performance of PHIF preconditioning for parabolic equations by
solving two examples of equation \eqref{eq1}: the heat equation and a logistic reaction-diffusion
equation. The PHIF preconditioning is compared with incomplete Cholesky preconditioning in terms of
the following quantities
\begin{itemize}
\item mem: the memory usage for PHIF and incomplete Cholesky factorizations;
\item $t_f$: the factorization time;
\item $t_s$: the average solve time for one time step, obtained from averaging over 100 time steps;
\item $n_i$: the number of CG iterations averaged over 100 time steps, with the relative residual
  equal to $10^{-12}$.
\end{itemize}

The only user-defined parameter in the PHIF factorization is the relative precision $\epsilon$ of the
interpolative decomposition, which is set to $10^{-3}$ and $10^{-6}$ in the numerical
experiments. Similarly, a drop tolerance $\epsilon$ is used for the incomplete Cholesky
factorization.

The MATLAB code for PHIF used for the numerical experiments is a modified
version of FLAM \cite{FLAM} to account for the block Jacobi preconditioning.

\paragraph{Example 1: the heat equation.} Consider first the 2d heat equation
\begin{align}
  \label{eq_heat}
  \frac{\partial u(x,t)}{\partial t} = \nabla\cdot(a(x)\nabla u(x,t)),\quad x\in\Omega=(0,1)^{2}
\end{align}
with the zero Dirichlet boundary condition and the initial condition equal to the sum of two
Gaussians,
\begin{align}
  \label{2d_h_u0}
  u_0(x,y) = e^{-((x-c_1)^2+(y-c_1)^2)/\sigma^2}+e^{-((x-c_2)^2+(y-c_2)^2)/\sigma^2},
\end{align}
where $c_1=0.35$, $c_1=0.65$ and $\sigma^2=0.05$. The diffusion coefficient field is set to
\begin{align}
  \label{2d_coeff}
  a(x,y) \sim \sum_{i=1}^{m}e^{-((x-x_i)^2+(y-y_i)^2)/\sigma_2^2},
\end{align}
rescaled and shifted to have values in the interval $[0.1,10]$, with $\sigma_2^2=0.005$,
$m=100$ and $x_i$ and $y_i$ being i.i.d. random variables sampled from the uniform distribution
$\mathcal{U}(0,1)$.

The time step size is set to $\Delta t= \Delta x = \frac{1}{N}$ and the numerical solution is
obtained for 100 successive time steps. Note that we want to set $\Delta t$ on the order of $\Delta
x$ to get a second order approximation. Since Crank-Nicolson is unconditionally stable we can select
a large time step, and since the initial condition is smooth, we don't observe numerical spurious
oscillations.

The numerical solution is evolved by solving \eqref{eq_CN2} at each time step using CG. The matrix
$A = I-\frac{\Delta t}{2}M$ is especially ill-conditioned for large time steps, therefore
preconditioning becomes necessary to reduce the number of CG iterations. Numerical results in Table
\ref{heat_2d} show a decrease on the number of CG iterations when using PHIF preconditioning instead
of incomplete Cholesky factorization, with the similar memory footprint.  The
  computation time per time step is approximately halved with the use of PHIF compared to incomplete
  Cholesky. Additionally, PHIF exhibits constant and problem size independent number of CG
  iterations, while the number of CG iterations for incomplete Cholesky scale with
  $O(N^{1/4})$. This results in almost linear $O(N)$ scaling of $t_s$ with PHIF preconditioning and
  $O(N^{1.25})$ scaling with incomplete Cholesky preconditioning.  PHIF also provides a good
approximation to the inverse for $\epsilon = 10^{-6}$, thus one could use the factorization to
directly solve the system of equations and bypass CG. For instance with
  $\epsilon=10^{-6}$ and $N=4095^2$, the PHIF factorization gives a solve error estimated as
  $\|I-\A\F^{-1}\|=7.9\times10^{-6}$ with randomized power iteration \cite{Dixon,Kuczynski} to
  $10^{-2}$ relative precision.

\begin{table}[htb]
\caption{Numerical results for the heat equation in 2d.}
\label{heat_2d}
\centering
\begin{tabular}{c|ccccc|ccccc}
\toprule
& \multicolumn{5}{|c|}{PHIF} & \multicolumn{5}{|c}{incomplete Cholesky}\\
 $N$ &$\epsilon$ & mem & $t_{f}$ & $t_s $ & $n_i$ & $\epsilon$ & mem & $t_f $ & $t_s $ & $n_i$\\
  & & (GB) & (s) & (s) & &  & (GB) & (s) & (s) & \\
\midrule
 $511^2$& & $0.19$ &\rev{$1.4e1$} &\rev{$1.7$}& $4.6$ &&$0.085$ &\rev{$1.7e{-1}$}&\rev{$1.6$}&$58.5$ \\
 $1023^2$ &$10^{-3}$&$0.80$&\rev{$6.2e1$}&\rev{$7.9$}&$5.2$&$10^{-3}$& $0.343$&\rev{$4.7e{-1}$}&\rev{$1.2e1$} &$97$\\
 $2047^2$& &$3.24$&\rev{$2.6e2$}&\rev{$3.6e1$}&$5.7$& &$1.37$ & \rev{$2.1$} & \rev{$7.4e1$} & $152.7$\\
 $4095^2$& &$13.0$&\rev{$9.1e2$}&\rev{$1.7e2$} & $5.8$& &$5.50$ &\rev{$7.7$}&\rev{$4.6e2$}& $226.8$\\
\midrule
  $511^2$& &$0.205$&\rev{$1.6e1$}&\rev{$8.6e{-1}$}&$2.3$ && $0.335$ &\rev{$1.4$}&\rev{$1.1$} & $11$ \\
 $1023^2$& $10^{-6}$&$0.834$&\rev{$6.6e1$}&\rev{$4.3$}&$2.7$& $10^{-5}$&$1.35$ &\rev{$6.1$}& \rev{$7.7$}& $16.4$\\
 $2047^2$& &$3.36$&\rev{$2.9e2$}&\rev{$2.1e1$}&$3$& &$5.4$ &\rev{$2.6e1$}& \rev{$4.3e1$}&$24.5$\\
 $4095^2$& &$13.5$&\rev{$9.8e2$}&\rev{$7.9e1$}&$3$& &$21.6$ &\rev{$1.1e2$}&\rev{$2.7e2$} & $38.2$\\
\bottomrule
\end{tabular}
\end{table}

Let us consider an analogous problem in three dimensions with $\Omega=(0,1)^3$ and the initial
condition equal to a Gaussian function,
\begin{align}
  \label{3d_u0}
  u_0(x,y,z) = e^{-((x-c)^2+(y-c)^2+(z-c)^2)/\sigma^2},
\end{align}
where $c=0.5$ and $\sigma^2=0.05$. The coefficient field is generated by
\begin{align}
  \label{3d_coeff}
  a(x,y,z) \sim \sum_{i=1}^{m}e^{-((x-x_i)^2+(y-y_i)^2+(z-z_i)^2)/\sigma_2^2}
\end{align}
rescaled and shifted to be within the interval $[0.05,20]$, with $\sigma_2^2=0.005$, $m=1000$
and $x_i$, $y_i$ and $z_i$ being i.i.d. random variables sampled from the uniform distribution
$\mathcal{U}(0,1)$. The time step is set to $\Delta t=0.1\Delta x$.

Numerical results are shown in Table \ref{heat_results3d}. Similarly to the 2d example, PHIF leads
to a reduction of CG iterations when compared to the threshold-based incomplete Cholesky. For
instance for $N=255^3$ with a tolerance $\epsilon = 10^{-6}$, CG with PHIF takes an average of 3
iterations, while CG with incomplete Cholesky takes 11.  While the factorization time
  of PHIF is high for this 3d example, the solve time per time step can be approximately halved with
  the use of PHIF as opposed to incomplete Cholesky.  For instance, the solve time per time step for
  PHIF is $1.8\times 10^2$ with $\epsilon = 10^{-6}$, while for incomplete Cholesky it is $3.8\times
  10^2$ for with $\epsilon = 10^{-4}$.  Additionally, experimentally PHIF exhibits constant number
  of CG iterations, while the number of CG iterations for incomplete Cholesky increases with the
  problem size $N$. This results in better scaling of $t_s$ with PHIF preconditioning than with
  incomplete Cholesky preconditioning, making PHIF better suited for large problem sizes.

\begin{table}[htb]
\caption{Numerical results for the heat equation in 3d.}
\label{heat_results3d}
\centering
\begin{tabular}{c|ccccc|ccccc}
\toprule
& \multicolumn{5}{|c|}{PHIF} & \multicolumn{5}{|c}{incomplete Cholesky}\\
 $N$ &$\epsilon$ & mem & $t_{f} $ & $t_s $ & $n_i$ & $\epsilon$ & mem & $t_f $ & $t_s $ & $n_i$\\
  & & (GB) & (s) & (s) & &  & (GB) & (s) & (s) & \\
\midrule
 $63^3$ & & $1.33$ &  \rev{$2.5e2$} & \rev{$3.8$} & $5.2$ & &$0.444$ & \rev{$6.0$} & \rev{$3.4$} & $14.5$\\
 $127^3$ & $10^{-3}$ & $14.1$ &  \rev{$3.2e3$} & \rev{$3.7e1$} & $5.9$ & $10^{-4}$ & $3.66$  &\rev{$1.9e1$} &\rev{$2.6e1$}&$22.9$ \\
 $255^3$& & $137$ & \rev{$1.3e4$} & \rev{$3.7e2$} & $6$ & &$29.7$ & \rev{$1.9e2$} & \rev{$3.8e2$} & $35$ \\
\midrule
  $63^3$ & &$1.78$ &  \rev{$8.7e1$} & \rev{$1.9$} & $3$ & &$3.54$ & \rev{$5.2e1$} & \rev{$6.0$} &  $5.8$\\
 $127^3$ & $10^{-6}$&$21.2$ &  \rev{$1.5e3$} & \rev{$2.0e1$} & $3$& $10^{-6}$ &$29.2$ & \rev{$8.2e2$} & \rev{$8.1e1$} & $8$\\
 $255^3$ & &$223$ &  \rev{$2.2e4$} & \rev{$1.8e2$} & $3$ & &$473$& \rev{$1.1e4$} & \rev{$1.4e3$} & $11$ \\
\bottomrule
\end{tabular}
\end{table}

\paragraph{Example 2: a logistic reaction-diffusion equation.} 
Consider now a 2d reaction-diffusion equation with logistic growth
\begin{align}
  \label{eq_nr}
  \frac{\partial u(x,t)}{\partial t} = \nabla\cdot(a(x)\nabla u(x,t)) + k_1u(x,t)
  \left(1-\frac{u(x,t)}{k_2}\right),\quad x\in\Omega=(0,1)^{2},
\end{align}
with $k_1=1$ and $k_2=10$, the zero Dirichlet boundary conditions and initial condition
\begin{align}
  \label{2d_u0}
  u_0(x,y) = \frac{2}{3\sqrt{2\pi\sigma^2}} e^{-((x-c)^2+(y-c)^2)/\sigma^2}
\end{align}
with $c=0.5$, $\sigma^2=0.05$. The coefficient field is set analogously to Example 1. This problem
is run with the same time step and the same number of time steps as in Example 1 and we observe no
numerical spurious oscillations in the numerical solution.

\begin{figure}[htb]
  \centering
  \begin{tabular}{cc}
    \includegraphics[width=0.45\textwidth,  clip]{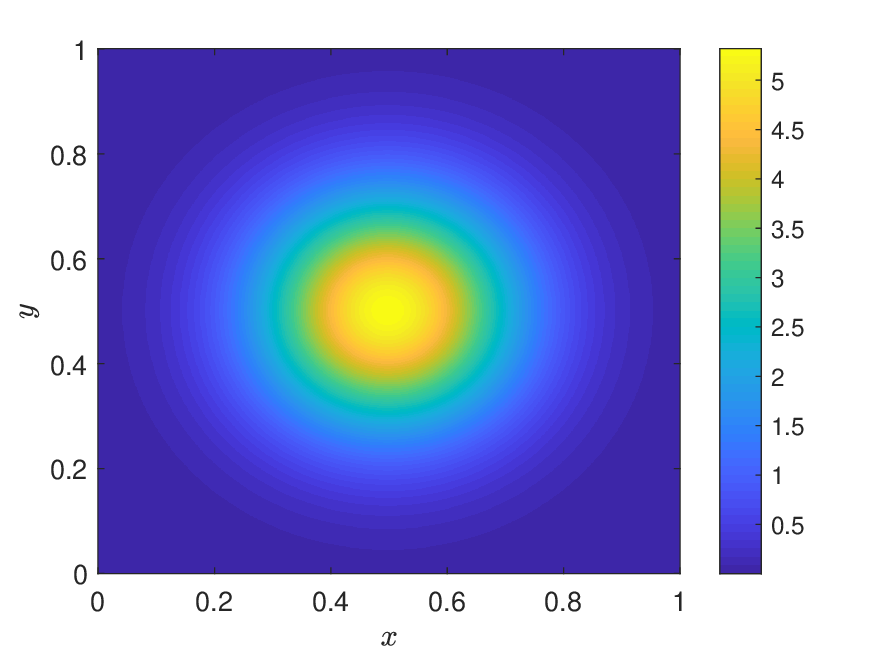} &
    \includegraphics[width=0.45\textwidth,  clip]{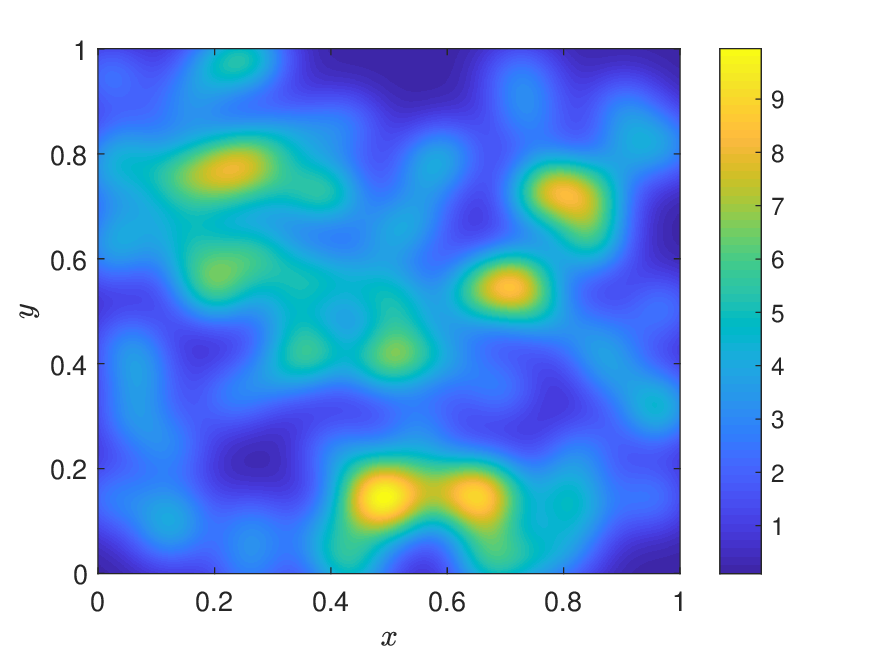}\\
    (a) $u_0(x)$ & (b) $a(x)$\\
    \includegraphics[width=0.45\textwidth,  clip]{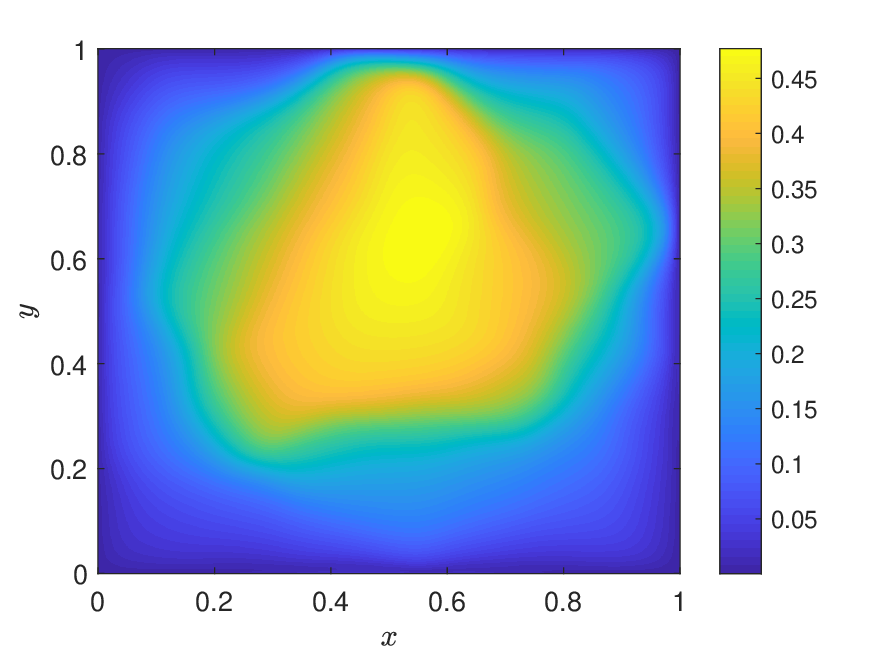} &
    \includegraphics[width=0.45\textwidth,  clip]{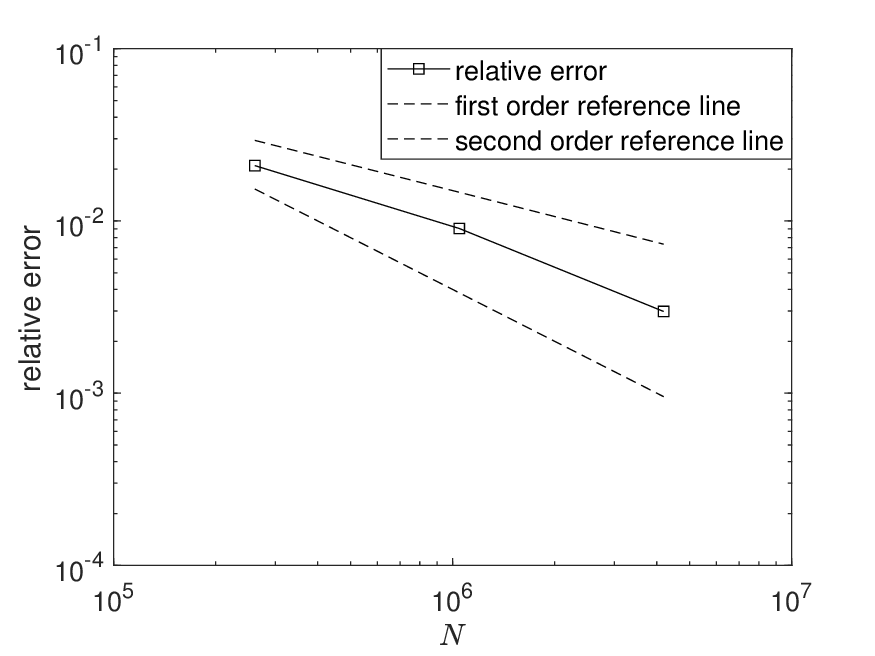}\\
    (c) $u(x)$ & (d) relative error
  \end{tabular}
  \caption{\label{fig:2d} Initial conditions (a), coefficient distribution (b), numerical solution
    (c) and relative error plot (d) for the 2d logistic reaction-diffusion equation.}
\end{figure}

The results are summarized in Table \ref{t_results}. Similarly to Example 1, we observe a decrease
on the number of CG iterations $n_i$ \rev{and on the computation time per time step $t_s$} when preconditioning with PHIF.  The initial condition and diffusion
coefficients are illustrated in Figure \ref{fig:2d} together with the numerical solution after 128
time steps for $N=4095^2$ and the relative error of the numerical solution $u^{(k=128)}$ for different problem sizes, which is close to
second order asymptotically.

\begin{table}[htb]
\caption{Numerical results for the 2d logistic reaction-diffusion equation}
\label{t_results}
\centering
\begin{tabular}{c|ccccc|ccccc}
\toprule
& \multicolumn{5}{|c|}{PHIF} & \multicolumn{5}{|c}{incomplete Cholesky}\\
 $N$ &$\epsilon$ & mem & $t_{f} $ & $t_s $ & $n_i$ & $\epsilon$ & mem & $t_f $ & $t_s$ & $n_i$\\
  & & (GB) & (s) & (s) & &  & (GB) & (s) & (s) & \\
\midrule
 $511^2$& & $0.19$ &\rev{$1.4e2$}&\rev{$2.3$}&$5.1$ &&$0.085$ &\rev{$1.5e{-1}$}&\rev{$1.8$}&$64.6$ \\
 $1023^2$&$10^{-3}$ &$0.81$&\rev{$5.7e1$}&\rev{$8.7$}&$5.6$& $10^{-3}$&$0.343$ &\rev{$5.7e{-1}$}&\rev{$1.3e1$}&$107.1$\\
 $2047^2$& &$3.24$&\rev{$2.3e2$}&\rev{$5.0e1$}&$6.5$& &$1.37$ & \rev{$2.4$} & \rev{$8.7e1$} & $167.1$\\
$4095^2$& &$13.0$&\rev{$8.5e2$}&\rev{$3.3e2$}&$6.3$& &$5.50$ & \rev{$7.8$} & \rev{$5.1e2$} & $252.9$\\
\midrule
 $511^2$& &$0.205$&\rev{$1.5$}&\rev{$9.2e{-1}$}&$2.6$ & &$0.335$ &\rev{$1.4$}&\rev{$1.3$}&$12.2$ \\
 $1023^2$ &$10^{-6}$ &$0.834$&\rev{$6.0e1$}&\rev{$4.3$}&$3$& $10^{-5}$ & $1.347$ &\rev{$6.2$}&\rev{$8.6$}& $18.5$\\
 $2047^2$& &$3.36$&\rev{$2.9e2$}&\rev{$2.8e1$}&$3.4$& &$5.39$ &\rev{$2.9e1$}&\rev{$5.7e1$}& $29.5$\\
 $4095^2$& &$13.5$&\rev{$9.8e2$}&\rev{$1.9e2$}&$4$& &$21.5$ & \rev{$1.2e2$} & \rev{$3.2e2$} & $40.8$\\
\bottomrule
\end{tabular}
\end{table}

For the three dimensional case with $\Omega=(0,1)^3$, we generate the initial condition $u_0(x)$ and
the coefficient field in the same way as the 3d case from example 1, with a multiplicative factor of
$(2\pi\sigma)^{-3/2}$ in \eqref{3d_u0} and $m=200$. We also set the time step to $\Delta t=0.1\Delta
x$ and the solution is evolved for 100 time steps.

\begin{figure}[h!]
  \centering
  \begin{tabular}{cc}
    \includegraphics[width=0.48\textwidth,  clip]{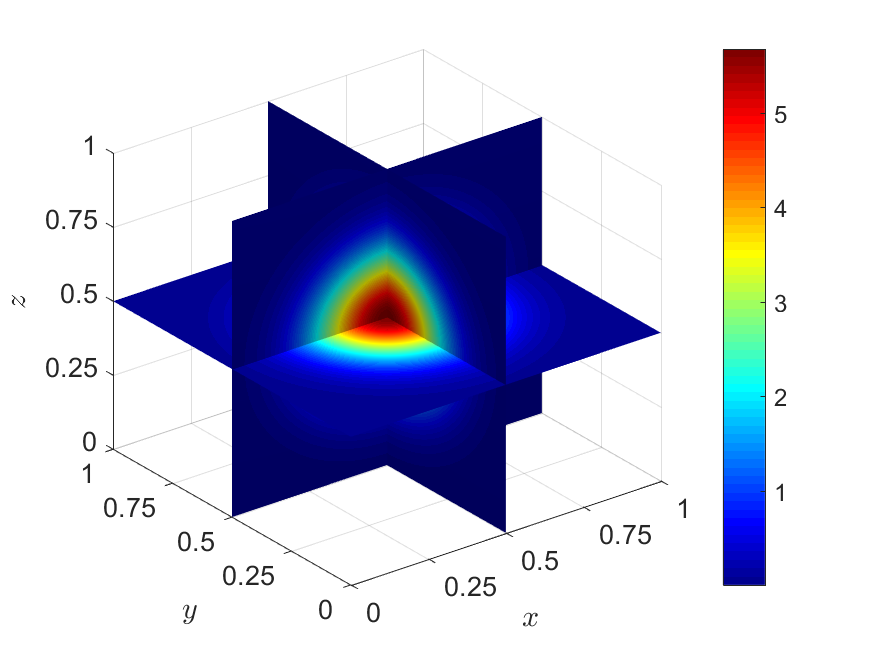} &
    \includegraphics[width=0.48\textwidth,  clip]{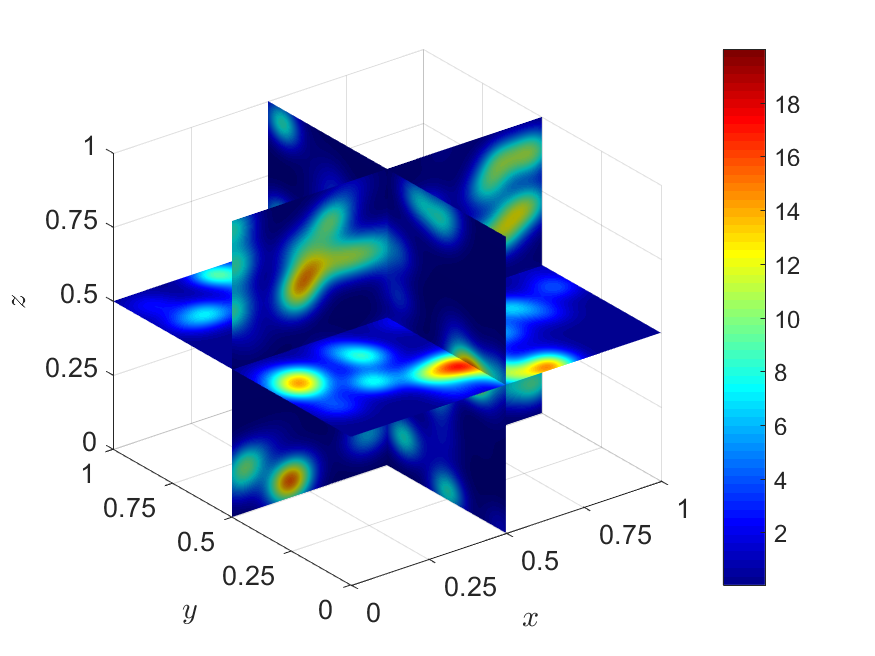}\\
    (a) $u_0(x)$ & (b) $a(x)$\\
    \includegraphics[width=0.48\textwidth,  clip]{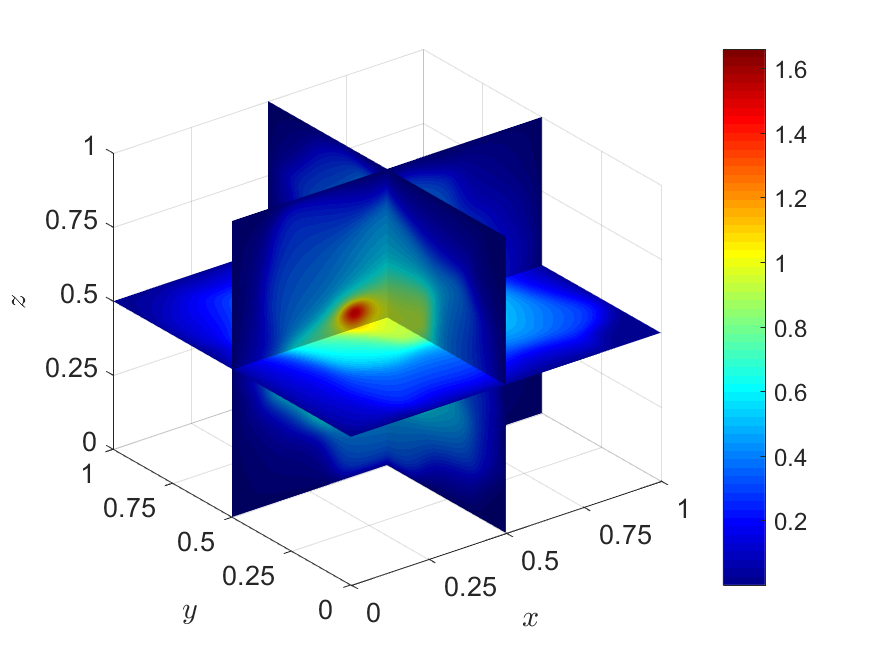} &
    \includegraphics[width=0.48\textwidth,  clip]{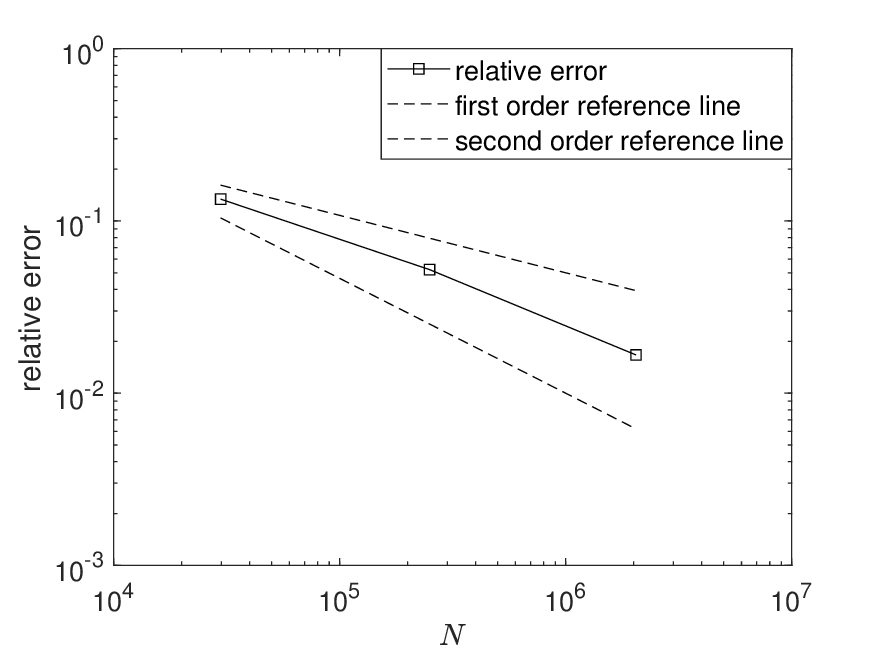}\\
    (c) $u(x)$ & (d) relative error    
  \end{tabular}
  \caption{\label{fig:3d} Initial conditions (a), coefficient distribution (b), numerical solution
    (c) and relative error plot (d) for the 3d logistic reaction-diffusion equation. }
\end{figure}

Numerical results are shown in Table \ref{t_results3d}. The initial conditions, coefficient field,
solution after 64 time steps for $N = 255^3$ and relative error of the numerical solution
$u^{(k=64)}$ for increasing problem sizes are depicted in Figure \ref{fig:3d}. We observe that the
error is close to second order asymptotically and CG converges with very few iterations using PHIF,
independently of $N$.

\begin{table}[htb]
\caption{Numerical results for 3d logistic reaction-diffusion equation}
\label{t_results3d}
\centering
\begin{tabular}{c|ccccc|ccccc}
\toprule
& \multicolumn{5}{|c|}{PHIF} & \multicolumn{5}{|c}{incomplete Cholesky}\\
 $N$ &$\epsilon$ & mem & $t_{f} $ & $t_s $ & $n_i$ & $\epsilon$ & mem  & $t_f $ & $t_s $ & $n_i$\\
  & & (GB) & (s) & (s) & &  & (GB) & (s) & (s) & \\
\midrule
 $63^3$ & & $1.31$ & \rev{$1.7e2$} & \rev{$4.1$} & $5$ & &$0.443$ &\rev{$2.4$} &\rev{$1.6$}& $14$ \\
 $127^3$ & $10^{-3}$ & $13.9$ &\rev{$9.6e2$} & \rev{$4.0e1$} & $6$ & $10^{-4}$ & $3.66$  &\rev{$1.6e1$} &\rev{$2.4e1$}&  $21$ \\
 $255^3$& & $135$ & \rev{$1.2e4$} & \rev{$3.7e2$} & $6$ & &$29.7$ & \rev{$1.5e2$} & \rev{$3.6e2$} & $32$ \\
\midrule
  $63^3$ & &$1.76$ & \rev{$1.8e2$} & \rev{$2.1$} & $3$ & &$2.77$ & \rev{$3.0e1$} & \rev{$4.2$} & $5.7$\\
 $127^3$ & $10^{-6}$&$21.2$ & \rev{$1.5e3$} & \rev{$2.3e1$} & $3$& $10^{-6}$ &$29.1$ & \rev{$5.2e2$} & \rev{$6.2e1$} &$7$\\
 $255^3$ & &$223$ & \rev{$2.3e4$} & \rev{$1.9e2$} & $3$ & &$237$& \rev{$7.6e3$} & \rev{$8.6e2$} & $10$ \\
\bottomrule
\end{tabular}
\end{table}

\section{Conclusions}\label{conclusions}

This note proposed an efficient preconditioner for solving linear and semi-linear parabolic
equations based on the hierarchical interpolative factorization in \cite{PHIF}. The preconditioned
CG iteration enjoys several advantages: (1) it provides a good approximate inverse that can be
applied very rapidly, (2) one only needs to construct HIF factorization once at the beginning, and
(3) applying the inverse approximation at each time step has linear cost. Computing the
factorization can be done in linear time as opposed to other more expensive factorizations such as
Cholesky or incomplete Cholesky. Well-suited for ill-conditioned matrices associated with large time
steps, the new preconditioner reduce the number of CG iterations significantly.

This approach can also be extended to solve other parabolic equations, for instance the
time-dependent fourth-order differential equations for studying the buckling plate or the clamping
plate problems in the plate theory \cite{Fishelov}. In such cases, HIF needs to use separators twice
as wide, when using the 9-point and 13-point finite differences stencils in 2d and 3d, respectively.

If the PDE has time-dependent coefficients or moving geometries, instead of
  constructing the PHIF factorization once at the beginning, one would need to compute the
  factorization at each time step. Such a computation can be expensive, especially for 3d
  geometries. However, if the coefficients change slowly with time, one can reuse the result by
  computing the factorization once every few time steps.

\section*{Acknowledgments}\label{acknowlegements}
The work of J.F. is partially supported by Stanford Graduate Fellowship. The work of L.Y. is
partially supported by U.S. Department of Energy, Office of Science, Office of Advanced Scientific
Computing Research, Scientific Discovery through Advanced Computing (SciDAC) program and the
National Science Foundation under award DMS-1818449.

\end{document}